\documentclass[12pt]{amsart} 
\usepackage{amsmath,amsthm, amssymb}
\usepackage{amsrefs}
\usepackage{times}
 \usepackage{hyperref}

\theoremstyle{plain}
\newtheorem*{main theorem}{Main Theorem}
\newtheorem*{theorem*}{Theorem}
\newtheorem*{emptytheorem*}{}
\newtheorem{theorem}{Theorem}

\newtheorem*{proposition*}{Proposition}

\newtheorem*{conjecture*}{Conjecture}
\theoremstyle{definition}
\newtheorem{definition}{Definition}
\newtheorem*{definition*}{Definition}
\theoremstyle{remark}

\newtheorem*{remark*}{Remark}

\newcommand{\dist}{\text{dist}}

\begin{document}
\title[Mixing in non-uniformly expanding maps]{Mixing and decay of correlations 
in non-uniformly expanding maps}

\author{Stefano Luzzatto}
\address{Department of Mathematics, 
Imperial College\\ London}
\email{\href{mailto:stefano.luzzatto@imperial.ac.uk}{stefano.luzzatto@imperial.ac.uk}}
\urladdr{\href{http://www.ma.ic.ac.uk/~luzzatto}
{http://www.ma.ic.ac.uk/\textasciitilde luzzatto}}
\date{Draft Version : September 2003} 

  \begin{abstract}
I discuss recent results on decay of correlations for non-uniformly 
expanding maps.  Throughout  the discussion, I address the 
question of why different dynamical 
systems have different rates of decay of correlations and how 
this may reflect underlying geometrical characteristics of the 
system. 
\end{abstract}
\thanks{Thanks to Gerhard Keller, Giulio Pianigiani, Mark Pollicott, Omri Sarig, and 
Benoit Saussol for 
useful comments on a preliminary version of this paper.}
\maketitle

\section{Introduction}

One of the basic questions of the theory of Dynamical Systems is to
describe the dynamics determined by the iterates of a map \( f: M\to M
\) on some space \( M \).  If \( x\in M \), we say that the (forward) orbit
of \( x \) is the set \( \theta (x) = \{f^{i}(x)\}_{i=0}^{\infty} \).
It is well known that even for  ``simple'' maps \( f \), the
topological and geometrical structure of orbits 
can be extremely complicated. However it is possible in some cases to 
obtain some remarkable results by focusing on the
statistical properties of these orbits rather than their precise 
geometrical and topological structure.   
In this survey we shall always assume that \( M \) is a smooth compact
Riemannian manifold of dimension \( d \geq 1 \). For 
simplicity we shall call the Riemannian volume \emph{Lebesgue 
measure} and denote it by \( m \) and assume that it is normalized so 
that \( m(M) = 1 \).   
Let \( f: M\to M \)  be a piecewise \( C^{2} \) map. 
For \( x\in M \) we let \( Df_{x} \) denote the derivative of \( f \) 
at \( x \) and define
\(
\|Df_{x}\| = \max\{\|Df_{x}(v)\|: v\in T_{x}M, \|v\|=1\}.
\)
We are  interested in formulating some 
\emph{expansion} conditions on  \( f \).  
The condition \( \|Df_{x}\|>1 \) implies that there is \emph{at
least} one vector which is expanded by \( Df_{x} \).  On the other
hand, the condition \( \|Df^{-1}_{x}\|<1 \), or equivalently \(
\|Df^{-1}_{x}\|^{-1}>1 \), means that \emph{all} vectors are contracted by the
inverse of \( Df_{x} \) and thus that \emph{all} vectors are
\emph{expanded} by \( Df_{x} \).  

\begin{definition} 
    We say that \( f \) is \emph{expanding on average}, or 
    \emph{non-uniformly expanding}
    \footnote{The \emph{non-uniformly expanding}
    here should be interpreted in the 
    sense of \emph{not necessarily} uniformly expanding
    rather than strictly \emph{not} 
    uniformly expanding so that the uniformly expanding case to be defined below 
    is a special case. The terminology is not optimal but 
    unfortunately the term \emph{expanding} is 
    traditionally interpreted to mean \emph{uniformly expanding}}
    if there exists \( \lambda>0 \) such that 
    \begin{equation*} \tag*{\( (*) \)}
    \liminf_{n\to\infty}\frac{1}{n}
    \sum_{i=0}^{n-1}\log \|Df^{-1}_{f^{i}(x)}\|^{-1}>\lambda
    \end{equation*}
    for almost every \( x\in M \).
\end{definition} 
An equivalent, and perhaps more intuitive, formulation is that for 
almost every \( x \) there exists a constant \( C_{x}>0 \) such that 
\[ 
\prod_{i=0}^{n-1} \|Df^{-1}_{f^{i}(x)}\|^{-1} \geq C_{x}e^{\lambda 
n}\quad \text{for every } n\geq 1.
\]
Thus every vector is expanded at a uniform 
exponential rate, although the constant \( C_{x} \) which can in 
principle be 
arbitrarily small, indicates that arbitrarily large number of 
iterates may be needed before this exponential growth becomes 
apparent. 
The expansivity is reflected at the level of the manifold by an 
exponentially fast divergence of nearby orbits, which produces the 
well known phenomenon of \emph{sensitive dependence on initial 
conditions}, in which 
small round-off errors increase very rapidly and orbits which start 
together  appear very quickly to have completely independent 
dynamics. This is often interpreted in a somewhat negative sense as 
leading to \emph{unpredictability} and \emph{chaos}. 
However it turns out that the very same 
expansivity property can give rise to some remarkable 
probabilistic structures which imply that there is a 
well-defined statistical coherence in the dynamics of different 
orbits. More precisely, it is possible to show that in many cases 
there exists an 
\emph{ergodic absolutely continuous invariant probability measure}. 
This implies that the asymptotic average distribution of orbits in 
space is the same for almost every initial condition \( x \). Thus 
the sensitive dependence on initial conditions implies that we cannot 
control the exact location of a given point at a given time, but the 
probabilistic structure implies that we can know that typical points 
will spend certain average amounts of time in certain regions of the 
space. 

The combination of the sensitive dependence on initial conditions and 
well defined average behaviour gives rise to the phenomenon of 
\emph{mixing}, which is a formalization of the notion that the 
dynamics, albeit completely deterministic, behaves to a large extent 
like a stochastic process. Leaving the precise definitions until 
later, we mention that the \emph{degree of stochasticity} can be 
quantified through the notion of 
 \emph{rate of mixing} or \emph{rate of decay of correlations} 
 which, in some sense, measure the rate at which the deterministic 
 system is approaching a stochastic process, or the speed at which 
 memory is lost.  It is known that a 
very wide range of rates can occur in different systems: 
from exponential to stretched exponential to polynomial to logarithmic. 
 However it is not completely clear \emph{why} different systems 
 exhibit different rates of decay of correlations. 
 The mathematics used to estimate the rate of decay of correlations
   for a particular system does not necessarily provide any intuitive 
   justification for why a particular rate occurs. 
   However some recent developments 
   are beginning to give us some insight into the mechanisms involved, 
   and point towards a subtle connection between the rate of 
   decay of correlations 
   and the finer geometrical structures associated to
   the dynamics. 
 The purpose of these notes is to 
 survey some recent results on rates of decay of correlations
 and to attempt, by bringing 
 these results together, to formulate some opinion on the general 
 question
\begin{quote}
\emph{What aspect of a non-uniformly expanding map determines the rate of 
mixing }?
\end{quote}

 For completeness we review the basic definitions related to invariant 
 measures and decay of correlations in Section \ref{basic notions}. 
 We then formulate the known results on decay for correlations for 
 uniformly expanding maps, maps with indifferent 
 fixed points, one-dimensional maps with critical points, and a class 
 of two dimensional Viana maps. These are pretty much all the 
 particular cases of maps for which results in this direction are 
 known.\footnote{We will concentrate in this survey on non-uniformly 
 expanding maps and not enter into any discussion of system which have 
 \emph{contracting} as well as \emph{expanding} directions, such as 
 \emph{partially} or \emph{non-uniformly} hyperbolic systems. The 
 spirit of the discussion remains valid in those examples 
 but there would be too much additional notation and different cases 
 to consider and the essential points of the discussion would be 
 overshadowed.}

This survey of known results suggests that one way of understanding 
the overall picture is the following. The exponential decay of 
correlations of uniformly expanding maps represents a kind of 
default rate, with slower rates occurring as a consequence of some 
intrinsic geometrical feature of a map which literally \emph{slows down} 
the process of orbits distributing themselves  
over the entire space. In some of the simpler examples this 
slowing down phenomenon is perfectly apparent and it is even possible 
to tweak some parameters to get essentially any desired 
(subexponential) rate. We do not yet have a general theory which 
confirms this point of view rigorously, although we shall present some 
recent results in this direction, which establish a 
 connection between the rate of decay of correlations and a certain 
 \emph{measure of non-uniformity} of the expansion, 
 related to the speed at which the 
\emph{liminf} in the definition of non-uniformly expanding maps is 
 attained.

\section{Basic notions}
\label{basic notions}
\subsection{Invariant measures}
The first key concept is that of an \emph{invariant probability measure}.
\begin{definition}
    A probability measure \( \mu \) on \( M \) is invariant under  \( f: M \to M \) if 
    \[ 
    \mu(f^{-1}(A)) = \mu(A)
    \]
    for every \( \mu- \) measurable set \( A \subset M \).
\end{definition}
Notice that this condition is not equivalent to the condition \(
\mu(f(A))=\mu(A) \) unless \( A \) is invertible.  In the general
non-invertible case we define 
\( 
f^{-1}(A) = \{x: f(x) \in A\}.
\)
By a general result of Ergodic Theory, some mild conditions on the map 
\( f: M \to M \) ensure that there exists at least one
invariant measure. However it is often the case, for example in the 
case of expanding maps which we are considering, that there are an
infinite number (for example every periodic orbit admits an invariant 
probability measure concentrated on that orbit). 
In this case there is an issue about which of
these measures is to be considered the most significant and what the
relation is between them. To discuss this notion in more depth we need 
some more definitions.

\subsection{Ergodicity and absolutely continuous measures.}

The relevance of invariant probability measure is that they capture 
the statistical features of the dynamics in a sense to be defined 
below. First we need one more definition.  
     
\begin{definition}
    Let \( \mu \) be an \( f \)- invariant measure.  We say that \( \mu
    \) is ergodic if there \emph{do not exist} measurable sets \( A, 
    A^{c} = M\setminus A  \) with 
    \( f^{-1}(A) = A \) (and thus also \( f^{-1}(A^{c})=A^{c} \)) and \( \mu(A)
\in (0, 1) \) (and thus also \( \mu(A^{c})= \mu(M\setminus A) \in 
(0,1) \)).
    \end{definition}
This means that if \( \mu \) is not ergodic, there are
(at least) two components of the ambient space \( M \) which never
interact; thus the dynamical properties of \( f \) on \( A \) 
completely independent of the dynamics \( A^{c} \).  
Any measure can be decomposed into a number of \emph{ergodic 
components} and, by  
a particular case of what is perhaps the most remarkable result of the 
Theory of Dynamical Systems, there is a strong statistical coherence 
in the dynamics of typical points associated to each ergodic component. 
    \begin{theorem*}[Birkhoff] 
    Let \( \mu \) be an ergodic invariant measure for \( f \).  
   Then, for any measurable set \( A\subset M \) and \( \mu \) almost 
   every \( x\in M \), we have
 \[ 
\frac{\#\{1\leq j \leq n : f^{j}(x) \in A\}}{n} \to \mu(A).
 \]
 \end{theorem*}
    Here \( \#\{1\leq j \leq n : f^{j}(x) \in A\} \) denotes the
    cardinality of the set of indices \( j \) for which \( f^{j}(x) \in A
    \).  The statement contains  two parts.  The first is the 
 convergence of the ratio on the left to some limit, a fact
    which is in itself quite remarkable.  If we interpret \( f^{j}(x) \in
    A \) as meaning that the event \( A \) occurs at time \( j \) for the
    initial condition \( x \) then this convergence means that there is a
    statistical pattern to the occurrence of such an event: as the number
    of iterates increases, the proportion of times for which the event occurs
    stabilises.  The second part of the statement is that this limit is
    precisely the \( \mu \)-measure of \( A \), sometimes stated in 
    the form  \emph{time averages
    equal space averages}. where ``space'' is measured according to the
    invariant measure \( \mu \).  
 
Of course, if the measure \( \mu \) is singular with respect to the
reference measure \( m \) on \( M \) given by the Riemannian volume, then this
result only concerns a set of zero volume and leaves us in the dark
about the dynamics of all other points. 
Things are better if the measure is \emph{absolutely
 continuous} with respect to the Riemannian volume.  
 \begin{definition}
     The measure \( \mu
 \) is absolutely continuous with respect to \( m \) if \( m(A) = 0 \)
 always implies \( \mu(A)=0 \). 
 \end{definition} 
 In particular, if \( \mu \) is absolutely continuous with respect to 
 \( m \), any set which has
 positive measure for \( \mu \) must also have positive measure for \( m \). 
Therefore, in this case, Birkhoff's Theorem immediately implies that there exists
 at least a positive volume set of points which are typical with
 respect to \( \mu \)\footnote{In many examples of
area-contracting diffeomorphisms any invariant measure is necessarily
supported on some \emph{attractor} which has zero volume, and
thus any invariant measure is necessarily singular with respect to the
Riemannian volume.  In many such cases there is a very sophisticated theory of
\emph{stable manifolds} which allows us to show that typical points in
\( M \) actually converge not just to the attractor as a whole but to
 typical points in the attractor, and are therefore typical with
 respect to the invariant measure; a conclusion which does not follow
 immediately from Birkhoff's Theorem.}.

 \subsection{Mixing and decay of correlations}  
An important first step in the analysis of non-uniformly 
 expanding maps in general, and in specific examples, is to show that 
 they admit an ergodic absolutely continuous invariant probability 
 measure. The next is to study the mixing properties of the map with 
 respect to that measure \cites{Hal44, Roh48, Ros56, Hol57, 
 KacKes58}. 
  
\begin{definition}
    An invariant probability measure \( \mu \) is mixing if 
    \[ 
    |\mu(A\cap f^{-n}(B))-\mu(A)\mu(B)| \to 0 
    \]
    as \( n\to \infty \), for all measurable sets \( A, B \subseteq M \).
    \end{definition}   
One way to interpret this condition is to notice that it 
 is  equivalent to the condition
  \[ 
  \left|\frac{\mu(A\cap f^{-n}(B))}{\mu(B)} -\mu(A)\right| \to 0
 \] 
as  \( n\to \infty \), for all measurable sets \( A, B \subseteq M \), with \( \mu(B) 
    \neq 0 \). 
Then one can think of
 \( f^{-n}(B) \) as a ``redistribution of mass'', notice that \( \mu(f^{-n}(B))
 = \mu(B) \) by the invariance of the measure.
 Then the mixing condition says that for large \( n \) the proportion
 of  \( f^{-n}(B) \) which intersects \( A \) is just proportional to the
 measure of \( A \).  In other words \( f^{-n}(B) \) is spreading
 itself uniformly with respect to the measure \( \mu \).  
A more probabilistic point of view is to think of \( \mu(A\cap
f^{-n}(B))/\mu(B) \) as the conditional probability of having \(
x \in A \) given the fact that \( f^{n}(x)\in B \), i.e. the
probability that the occurrence of the event \( B \) today is a
consequence of the occurrence of the event \( A \) exactly \( n \)
steps in the past.  The mixing condition then says that this 
probability converges precisely to
 the probability of \( A \). Thus, asymptotically,
there is no causal relation between the two events.
 
 The mixing condition can be written in an
 integral form as
 \[ 
 \left|\int \chi_{A\cap f^{-n}(B)} d\mu - \int \chi_{A} d\mu \int
 \chi_{B} d\mu\right| \to 0
 \]
 or even 
 \[ 
  \left|\int \chi_{A} (\chi_{B}\circ f^{n}) d\mu - \int \chi_{A} d\mu \int
  \chi_{B} d\mu\right| \to 0
 \]
 This last formulation now admits a natural generalization by replacing
 the characteristic functions with arbitrary measurable functions. 
 \begin{definition}
     For real valued measurable functions \( \varphi, \psi : M \to \mathbb
     R \) we define the correlation function 
     \[ 
     \mathcal C_{n}(\varphi, \psi) = \left|\int \psi (\varphi\circ f^{n})
     d\mu - \int \psi d\mu \int \varphi d\mu\right|
     \]
 \end{definition}
 In this context, the functions \( \varphi \) and \( \psi \) are often
 called \emph{observables}.  Assuming that the measure \( \mu \) is mixing, a
 natural question is whether \( \mathcal C_{n}(\varphi, \psi) \) tends 
 to \( 0 \) for general observables and, if so,
 whether there is any particular \emph{speed} at which the correlation
 function decays, i.e. does there exist a sequence \( \{\gamma_{n}\} \) 
 with \( \gamma_{n}\to 0 \) as \( n\to \infty \) (for example \(
 \gamma_{n}=e^{-\alpha n} \) or \( \gamma_{n}=n^{-\alpha} \)),
 depending only on the map \( f \) and on the class of admissible
 observables, such that for any two admissible \( \varphi, \psi \)
 there exists a constant \( C_{\varphi, \psi} \) such that
 \[ 
 \mathcal C_{n}(\varphi, \psi) \leq C_{\varphi, \psi} \gamma_{n}
 \quad \forall \ n.
 \]
 It turns out that the answer is negative if the class of admissible
 observables is too large, e.g. \( L^{\infty}(\mu) \) which in
 particular contains the characteristic functions.  Indeed, it is a
 classical result that for any given sequence \( \{\gamma_{n}\} \) it
 is possible to choose measurable sets \( A, B \) such that the
 correlation function \( \mathcal C_{n}(\chi_{A}, \chi_{B}) \) decays
 slower than the rate determined by the sequence \( \{\gamma_{n}\} \). 
 However, if the class of
 admissible observables is restricted to functions with some
 regularity, e.g. continuous functions with some conditions on the
 modulus of continuity such as H\"older continuous functions, it turns out
 that for a large and significant class of dynamical systems it is
 indeed possible to speak of a particular \emph{rate of mixing} or
 \emph{rate of decay of correlations}.  Unless explicitly mentioned 
 below, we shall always assume that we are dealing with H\"older 
 continuous observables.

 \section{Uniformly Expanding Maps}

 \subsection{The smooth case}
The definition of (non-uniformly) expanding
includes as a special case the classical uniformly expanding case.
 \begin{definition}
     We say that \( f \) is uniformly expanding if there exist constants \(
     C, \lambda > 0 \) such that for all \( x\in M \), all \( v\in T_{x}M
     \), and all \( n\geq 0 \), we have 
     \[ 
     \|Df^{n}_{x}(v)\|\geq Ce^{\lambda n}\|v\|.
     \]
     \end{definition}

In this case, as in most of the cases which will be discussed below, 
the proof of the existence of an absolutely continuous invariant 
measure historically predates the estimates for the rate of decay of 
correlations, and can usually be obtained through significantly 
simpler arguments. 

\begin{theorem}[\cites{Ren57, Par60, Gel59, Rue68, Ave68, 
    KrzSzl69, Wat70, Las73}] 
 Let \( f: M \to M \) be \( C^{2} \) uniformly expanding. Then there
     exists a unique absolutely continuous mixing 
     invariant probability measure \( \mu
     \). 
\end{theorem}

\begin{theorem}[\cites{Sin72, Per74, Bow75, Rue76}]
 Let \( f: M \to M \) be \( C^{2} \) uniformly expanding.
 Then the correlation function decays exponentially fast.
 \end{theorem}
 
 Several of the arguments used in the proofs of these results 
 rely on the fact that smooth 
 uniformly expanding systems have a geometrical \emph{Markov} 
 structure, that is there 
 exists a   
 partition of \( M \) modulo sets of measure 
 0 (finite in this particular case), such that \( f \) is \( C^{2} \) 
 on each partition element and maps each element 
 to the whole manifold, or to some suitably large union 
 of other partition elements. The results then generalize quite 
 naturally to the more general Markov case, even if the partition is 
 countable, as long as some mild technical conditions are satisfied
 \cites{You98, Bre99, You99, Mau01a,  Mau01b}. 
 
 \subsection{The piecewise smooth case}
 
The general (non-Markov) piecewise expanding case is significantly 
more complicated and even the existence of an absolutely continuous 
invariant measure is no longer guaranteed \cites{LasYor73, 
GorSch89, Tsu00a, Buz01a}. The main problems 
lie in the fact that the images of the discontinuity set 
can be very badly distributed and 
cause havoc with any kind of \emph{structure}. In the Markov case this 
does not happen because the set of discontinuities gets mapped to 
itself by definition.  Also the possibility 
of components being \emph{translated} in different directions can 
destroy on a global level 
the local expansiveness given by the derivative. 
Moreover, where results exist for rates of decay of correlations, they do 
not always apply to the case of H\"older continuous observables, as 
technical reasons sometimes require that different functions spaces be 
considered which are more compatible with the discontinuous nature of 
the maps. We shall not explicitly comment on the particular classes 
of observables considered in each case. 

In the one-dimensional case these problems are somewhat more 
controllable and relatively simple conditions guaranteeing the 
existence of an ergodic invariant probability measure can be 
formulated even in the case of a countable number of domains of 
smoothness of the map. These essentially require that the size of the 
image of all domains on which the map is \( C^{2} \) strictly 
positive  and that certain conditions on 
the second derivative are satisfied
\cites{LasYor73, Adl73, Bow77, Bow79}. 
 In the higher dimensional case, the situation is considerably more 
 complicated and there are a variety of possible conditions which can 
 be assumed on the discontinuities. The conditions of 
 \cite{LasYor73} were generalized to the two-dimensional context in 
 \cite{Kel79} and then to arbitrary dimensions in \cites{GorBoy89, 
 Tsu01b, Buz00a}. There are also several other papers which 
 prove similar results under various conditions, we mention 
 \cites{Alv00, BuzKel01, BuzPacSch01, Buz01c, Sau00, BuzSar}.
 In \cites{Buz99a, Cow02} it is shown that conditions 
 sufficient for the existence of a measure are \emph{generic} in 
 a certain sense within the class of piecewise expanding maps. 
 
Exponential decay of correlations has also been proved for non-Markov 
piecewise smooth maps, although again the techniques have had to be 
considerably generalized. In terms of setting up the basic arguments 
and techniques, a similar role to that played by 
\cite{LasYor73} for the existence of absolutely continuous invariant 
measures might be attributed to \cites{Kel80, HofKel82, Ryc83} 
for the problem of decay of correlations in the one-dimensional 
context. More recently, alternative approaches have been proposed and 
implemented in \cites{Liv95a, You98}. The approach of \cite{You98} 
has proved particularly suitable for handling some higher dimensional 
cases such as 
\cite{BuzMau03} in which assumptions on the discontinuity set 
are formulated in terms of the \emph{topological pressure} of this 
set and exponential decay of correlations is proved,
and \cites{AlvLuzPin, AlvLuzPin03} in which 
the assumptions on the discontinuity set 
are formulated  as \emph{geometrical non-degeneracy} 
assumptions, which are essentially conditions on the first and second 
derivatives of the map near the discontinuity, and \emph{dynamical } 
assumptions on the \emph{rate of recurrence} 
of typical points to the discontinuities. 
We remark however that the estimates in 
\cites{AlvLuzPin, AlvLuzPin03} are sub-exponential (and sub-optimal) 
and clearly weaker than the exponential decay of correlations 
obtained in
\cite{BuzMau03}, 
 perhaps because they are carried out in a general 
framework which allows the
possible presence of critical points (this will be discussed in more 
detail in  Section \ref{general theory}). It would be interesting to 
know how the assumptions of \cites{AlvLuzPin, AlvLuzPin03} and 
those of \cite{BuzMau03} are related.

 \section{Almost uniformly expanding maps}
  \label{almost expanding}
 \subsection{Neutral fixed points}
 Perhaps the simplest way to weaken the uniform expansivity condition 
 is to consider a one-dimensional map which is expanding, i.e. \( 
 |f(x)|>1 \), everywhere 
 except at some fixed point \( p \) at which \( f'(p) =1 \). 
The fixed point \( p \) is still repelling, but nearby points remain close to \( p \)
much longer than they would if the derivative were \( >1 \). 
On the other hand the dynamics away from the 
 fixed point is uniformly expanding and orbits tend to distribute themselves
over the whole interval quite quickly. 
 Thus the overall 
 effect is that orbits tend to spend a long time \emph{trapped} 
 in a neighbourhood of the fixed point with relatively short 
 \emph{bursts of chaotic activity} outside this neighbourhood. 
 This is characteristic of the phenomenon of \emph{intermittency} 
 which appears to be  common in many natural phenomena and is 
an important feature for example in the theory of \emph{Self-Organized 
 Criticalities} \cites{Bak97, Jen98}. Indeed, this was one of the 
 motivations which led to the study of maps with this kind of 
 characteristics, see \cite{ManPom80}.
  
 More formally, 
 we suppose that there  exists a  partition 
 \( \mathcal P \) of \( [0,1] \) into a finite number of subintervals 
 and that
     \( f \) is \( C^{2} \) in the interior of each partition element 
     with a \( C^{1} \) extension to the boundaries and that the derivative is 
     strictly greater than 1 everywhere except at a 
  fixed point \( p \) (which for simplicity we can assume lies at the origin) 
  where \( f'(p)=1 \). For the moment 
  we assume also a strong Markov property: 
  each partition element is mapped bijectively to the whole interval. 
 First of all we want to focus on the consequences of 
  the presence of the \emph{neutral} (or \emph{indifferent)} fixed point \( p \). 
  For definiteness, let us suppose that 
 on a small neighbourhood of 0 the map takes the form 
  \[ 
  f(x) \approx x+x^{2}\phi(x)
  \]
  where \( \approx \) means that the terms on the two sides of the 
  expression as well as their first and second order derivatives 
converge as \( x\to 0 \). We assume moreover that
 \( \phi \) is \( C^{\infty} \) for \( x\neq 0 \); 
  the precise form of \( \phi  \) determines the precise degree of \emph{neutrality} 
  of the fixed point, 
and in particular  affects the second derivative \( f'' \). It turn out that it 
plays a crucial role in determining the mixing properties and even 
the very existence of an absolutely continuous invariant measure. 

\subsection{Loss of mixing}

The following result shows that the situation is drastically different 
from the uniformly expanding case.

\begin{theorem}\cite{Pia80}
    If \( f \) is \( C^{2} \) at the neutral fixed point,  
    then \( f \) does 
    not admit any absolutely continuous invariant probability measures.
 \end{theorem}

 This case occurs if, for example, \( \phi(x) \equiv 1 \). 
It is interesting to note that 
\( f \) has the same topological behaviour as a uniformly 
expanding map,  
 typical orbits continue to wander densely on the whole interval, but 
 the proportion of time which they spend in various regions tends to 
 concentrate on the fixed point, so that, asymptotically, typical 
 orbits spend all their time near \( 0 \). 
 It turns out that in this situation there 
 exists an infinite (\( \sigma \)-finite) 
 absolutely continuous invariant measure which 
 gives finite mass to any 
 measurable set not containing the fixed point 0 and infinite 
 mass of any neighbourhood of \( p  \)
 \cite{Tha83}. 
 
 \subsection{Non-uniform expansivity}
The situation  changes if 
 we relax the condition that \( f \) be \( C^{2} \) at \( p \) and 
 allow the second derivative \( f''(x) \) to diverge to 
 infinity as \( x\to p \). This means that the derivative 
 increases 
  quickly near \( p \) and thus nearby points are repelled at a 
  faster rate. Although this is apparently a very 
 subtle change, it makes all the difference. With some mild 
 conditions on the rate of divergence of \( f'' \) near \( p \) it is 
 possible to recover the existence of an absolutely continuous 
 probability measure \( \mu \) \cites{Pia80, HuYou95}. 
 Typical points still spend a large proportion 
 of time near \( p \) but they now also spend a positive proportion of 
 time in the remaining part of the space.  In particular they are 
 non-uniformly expanding:  by a simple
 application of Birkhoff's Ergodic Theorem to the function \( \log |f'(x)| \), 
we have that, for \( \mu \)-almost every \( x \), 
 \[ 
 \lim_{n\to\infty}\frac{1}{n}\sum_{i=0}^{n-1}\log |f'(f^{i}(x))| \to \int 
 \log|f'| d\mu > 0.
 \]
  The fact that \( \int \log
 |f'|d\mu > 0 \) follows from the simple observation that \( \mu \) is
 absolutely continuous and finite,  and that \( \log |f'| > 0 \) 
 except at the neutral fixed point. 
 
 \subsection{Polynomial decay of correlations}
If \( \phi \) is  of the form 
 \[ \phi(x) = x^{-\alpha}\quad\text{for some } \alpha\in (0,1), \]
then we have the following 
 \begin{theorem}\label{neutral polynomial}
     \cites{Hu01, Iso99, LSV98, You99, Sar02}
     \( f \) admits a mixing (in particular ergodic) invariant probability 
     measures with decay of correlations
     \[ 
     \mathcal C_{n} = \mathcal O(n^{1-\frac{1}{\alpha}}).
     \]
     \end{theorem}
     Thus, the existence of an absolutely continuous invariant measure as in the 
     uniformly expanding case 
     has been recovered, but the exponential rate of decay of correlation has not. 
Thinking of the rate of decay of correlation as related to the mixing 
process, and in particular to the speed at which 
 the mass gets redistributed by the dynamics, we can think of 
the indifferent fixed point as having the effect of 
 \emph{slowing down} 
 this process by trapping nearby points for disproportionately long time.   
     In fact we have a \emph{slower rate of decay} if \( \alpha \) is larger,  
 which corresponds to the second derivative \( f'' \) diverging
 more slowly as the fixed point is approached.  Geometrically, 
 the order of the tangency between the graph of \( f \) and 
 the diagonal increases with \( \alpha \), 
 thus the fixed point is becoming \emph{less repelling}, points tend to 
   remain trapped for longer, and the rate of 
   decay of correlations is slower.   

   We remark that the estimate in \cite{Sar02} shows that the rate 
   stated is optimal by obtaining lower bounds as well as upper 
   bounds. Most of the results stated in this survey are upper 
   bounds. A brief discussion of this issue is given at the end of 
   the paper. 
   
   \subsection{Logarithmic decay of correlations}
This intuitive explanation for the connection between 
the order of the fixed point and the rate of decay of correlations
is supported by recent work \cite{Hol03} in which the situation is taken to an extreme 
by considering very general functions \( \phi \) satisfying a \emph{slowly varying condition} 
\cite{Aar97}. This gives a range of possible rates of decay of correlation, 
including \emph{logarithmic, intermediate logarithmic} and \emph{intermediate polynomial}. 
As an example, if we let \( \phi \) be of the form:
 \[ 
 \phi(1/x) = \log x \log^{(2)}x\ldots\log^{(r-1)}x (\log^{(r)}x)^{1+\alpha}
 \]
 for some \( r\geq 1,\) \(\alpha\in (-1,\infty) \) 
 where \( \log^{(r)}=\log\log\ldots\log \) 
 repeated \( r \) times, we get the following result: 

 \begin{theorem}\cite{Hol03}\label{neutral logarithmic}
\( f \) admits a mixing absolutely continuous invariant probability measure with 
decay of correlations
\[ 
\mathcal C_{n}=\mathcal O (\log^{(r)}n)^{-\alpha}.
\]
\end{theorem}

We remark that the methods of  \cite{You99} for the proof of Theorem 
\ref{neutral polynomial} apply to a 
considerably larger class of maps than those explicitly defined here. 
In particular we can consider 
maps with any arbitrary finite number of neutral fixed or periodic orbits 
and, most importantly, 
the Markov condition can be relaxed by adding some mild additional conditions 
on the expansivity 
(i.e. assuming that \( f'\geq\mu>2 \) on the partition elements which do 
not contain the neutral fixed point). 
The arguments of \cite{Hol03} are based on generalizations of the methods of 
\cite{You99} and are therefore very likely to extend to 
give logarithmic decay of correlations in these additional cases as well. 
Recently there have been some generalizations of the results above to 
higher-dimensional situations, see
\cites{PolYur01a, Hu01}.

 \section{One-dimensional maps with critical points}
 
 \subsection{Unimodal and multimodal maps}
 We now consider another class of systems which can also exhibit
 various rates of decay of correlations, but where the mechanism for
 producing these different rates is significantly more subtle.  We 
 consider
 the class of \( C^{3} \) maps \( f: [0,1] \to [0,1]\) with some
 finite number of non-flat critical points.  We say that \( c\in [0,1]
 \) is a critical point if \( f'(c) = 0 \); the critical
 point is non-flat if there exists an \( 0< \ell<\infty \) called the
 \emph{order} of the critical point, such that \( |f(x)| \approx
 |x-c|^{\ell} \) for \( x \) near \( c \);  \( f \) is
 unimodal if it has only one critical point, and multimodal if it has
 more than one. We shall allow \( f \) to have an arbitrary finite number of
 critical points but assume that they all
 have the same order.  We shall also assume a standard technical
 condition called \emph{negative Schwarzian derivative} which is a
 kind of convexity assumption on the derivative of \( f \). 
Although the results we shall describe apply in significant generality,
they are already interesting and highly non-trivial 
in the case of maps belonging to the
well known logistic family
 \[ 
 f_{a}(x) = ax(1-x)
 \]
which has a unique critical point of order \( \ell=2 \) (and 
satisfies the negative Schwarzian derivative condition). 

\subsection{Non-uniform expansion}
As mentioned
above, the questions concerning the existence and ergodicity of an
absolutely continuous invariant probability measure are of a 
somewhat different
nature from those which concern the statistical properties of such
measures.  A first
important observation is that the expanding properties of the maps,
which play an important role in the analysis of the two classes of
examples mentioned above, are not at all obvious here.   
Around the critical point there is a region in which the
derivative is arbitrarily small, and close to the boundaries of the
interval \( [0,1] \) there are regions in which the derivative is
relatively large.  The non-uniform expansivity condition is still
possible in principle but depends on typical orbits spending on
average more time in the expanding region than in the contracting
region near the critical points.  The first observation that this
actually does happen goes back to Ulam and von Neumann \cite{UlaNeu47} in
which they show that the so called \emph{top} quadratic map, \( f(x) =
4x(1-x) \) (they actually use a different but equivalent
representation of the quadratic family), is \( C^{1} \) conjugate to a
uniformly expanding map.  This allows them to show explicitly that
there exists an absolutely continuous invariant probability measure
and to show that the map is therefore non-uniformly expanding in the
sense of our definition.

The approach of Ulam and von Neumann however 
relies on some special characteristics of the particular map they
consider. It can be generalized to work for some relatively small 
(albeit infinite) set of parameter values, although this 
generalization 
depends on results from the theory of one-dimensional dynamics 
which were not available until the late 70's and the beginning of the  
80's, 
many decades after the work of 
\cite{UlaNeu47}. 
This theory has spawned a huge amount of research and
it falls well beyond the scope of this survey to enter into the
technical details of the differences between one result and another. 
We really just mention the basic principle of what has been
understood, which is that a substantial amount of information about
the overall expansivity properties of the map is contained in the
expansivity properties of the orbits of the critical points, i.e. in
the behaviour of the sequence \( D_{n}(c) = |(f^{n})'(f(c))|\).
\footnote{The fact that the derivative is calculated in the critical 
value \( f(c) \) rather than the critical point \( c \) often causes 
confusion to people not used to working with maps with critical 
points. It is sufficient to observe however that by the chain rule we 
have \( (f^{n})'(x) = f'(x) f'(f(x))\ldots(f'(f^{n-1}(x))\) and 
therefore this would always be equal to zero if we choose \( x=c \) 
since \( f'(c)=0 \) by definition. Calculating the derivatives at \( 
f(c) \) on the other hand gives an accurate reflection of the behaviour 
of the critical orbit.} 
In the unimodal case, there are several papers
proving the existence of absolutely continuous invariant probability 
measures 
under weaker and weaker assumptions: finite critical orbit 
\cite{Rue77}, non-recurrent critical point 
 \cite{Mis81} (both of these conditions imply 
 exponential growth of \( D_{n} \)), 
exponential growth \cites{ColEck83, NowStr88}, a
summability condition \cite{NowStr91}, and
more recently the remarkable paper \cite{BruStrShe03} which supersedes 
all previous results 
by showing that an absolutely continuous
invariant probability measure exists under the simple condition that
\( D_{n}\to \infty \) without any further assumption about the rate
of growth of this derivative.  The multimodal case has proved
significantly harder and was first addressed in the following
\begin{theorem}\cite{BruLuzStr03}
\label{existence}
If \( f \) satisfies 
\begin{equation*}
\sum_n D_n^{-1/(2\ell-1)} < \infty \
\end{equation*}
for each critical point \( c \), then there exists an \( f \)-invariant
probability measure \( \mu \) absolutely continuous with respect to
Lebesgue measure.
\end{theorem}
The assumptions have since been weakened to the
summability condition \( \sum D_{n}^{-1/\ell}<\infty \) and to
allow the possibility of critical points of different orders\cite{BruStr01}. 

It is worth mentioning that, unlike the examples of the previous 
sections, the conditions on the growth of the derivative along the 
critical orbit which define this class of examples, 
are generally not directly verifiable since they involve the full 
forward orbit of the critical point. Moreover, these conditions 
are extremely unstable : 
arbitrarily small changes in the 
parameter can destroy the delicate balance between the number of 
iterates spent in the contracting region and in the expanding regions 
which cause the orbit to exhibit derivative growth on 
average.\footnote{Strictly speaking this has only been 
proved in certain unimodal cases \cites{GraSwi97, Lyu97, Koz03}, but 
it is widely believed to be a very general fact.} 
Non-trivial arguments are thus required to show that the 
appropriate growth conditions are satisfied for a significant set of 
parameter values. It turns out that in generic one-parameter families 
there exist \emph{nowhere dense} sets of \emph{positive Lebesgue measure}
\footnote{Thus the non-uniform expansion property in the context of 
one-dimensional maps with critical points occurs for a set of 
parameters which is
\emph{topologically small} but \emph{measure-theoretically large}.}
for 
which the corresponding maps 
admit an absolutely continuous invariant probability measure 
\cites{Jak81, BenCar85, Ryc88, ThiTreYou94, MelStr93, Tsu93a, Tsu93, 
Luz00}. There are also some generalizations to families of piecewise 
smooth maps with critical points \cites{LuzVia00, LuzTuc99}.

\subsection{Basic strategy} 
These results rely
in a crucial way on the (one could say \emph{miraculous}\ ) 
property of such one-dimensional
maps that, as long as there are no periodic attractors, the
dynamics outside any fixed neighbourhood of the critical set is
uniformly expanding \cites{Man85, Man87}, and on the remarkable 
observation that if the derivative along the critical orbit is bounded 
away from 0 there can be no attracting periodic orbits (as an 
immediate consequence of the fact that any attracting periodic orbit 
must have a critical point in its basin of attraction, see \cite{ 
MelStr93}).  In general this 
expansivity will degenerate as the size of the neighbourhood 
tends to zero, but the crucial fact remains that we
can divide the interval into some arbitrarily small neighbourhood \(
\Delta \) and its complement \( [0,1]\setminus \Delta \); since the
dynamics is uniformly expanding on \( [0,1]\setminus\Delta \) it is
necessary, and to some extent sufficient, to control the dynamics and 
recurrence of points in \( \Delta \).  Orbits that pass through
\( \Delta \) may
lose a lot of the expansion which they had previously
accumulated outside \( \Delta \), but if they do not fall in \( \Delta
\) too often, and if, when they do, they do not fall too close to the
critical point, it may happen that this loss is not sufficient to
destroy the expansivity completely.  One possible strategy for
controlling the extent of this loss of expansion is to take advantage
of the fact that points in \( \Delta \) 
are mapped extremely close to a critical value and
thus inevitably \emph{shadow} the orbit of the critical value for a
certain amount of time.  The closer the point is to the critical point
\( c \), the longer the shadowing time.  During this time, the point
has essentially the same behaviour as the critical value and in
particular has the same pattern of derivative growth.  Thus if we
know that the critical value has some expansivity
properties and we can show that our original orbit shadows it for long
enough, we can conclude that by the end of this shadowing period it
will have regained the expansion it lost by having an iterate in \(
\Delta \).
 
The details of this argument rely on a balance between the rate of 
growth of the derivative along the critical orbits and the number of 
iterates for which a point \( x\in \Delta \) shadows this orbit 
which in turn depend on how close \( x \) is to the critical point 
which in turn determines how much expansion is lost at the first 
iterate. 
If the derivative 
is growing very fast, e.g. exponentially, then \( x \) tends to 
get \emph{pushed away} much faster, and the shadowing period is 
relatively short. However, since the derivative is growing 
very fast, this length of time is sufficient to recover the loss of 
derivative incurred at the return in \( \Delta \). If the derivative 
along the critical orbit is growing more slowly, e.g. at a polynomial 
rate, then nearby points tend to get pushed away more slowly and the 
shadowing lasts significantly longer. This is good from the point of 
view of recovering expansion as more time is needed due to the very 
same slow derivative growth along the critical orbit. As far as the 
existence of of an absolutely continuous invariant measure is 
concerned, the estimates do indeed work, and we can show that even 
for a relatively slow rate of growth of the derivative along the 
critical orbits, such as that given by the summability condition in 
the Theorem, the shadowing of the critical orbit is sufficient to 
compensate the small derivative at returns to \( \Delta \) for almost 
all orbits, giving rise to a map which satisfies the non-uniform 
expansivity condition \cite{BruLuzStr03}. 

\subsection{Rate of mixing}
This approach becomes particularly interesting when 
addressing the question of the decay of correlations. It turns out 
that the \emph{existence} of critical points is much less of an issue than 
the particular \emph{rate of growth} of the derivative along the critical 
orbits, and that the conceptual picture is much more similar to the 
case of expanding maps with indifferent fixed points discussed above, 
than would appear at first sight.  Indeed, we can think of the case 
in which the rate of growth of \( D_{n} \) is subexponential as a 
situation in which the critical orbit is in some sense \emph{neutral} 
or \emph{indifferent} and in this respect very similar to the neutral 
or indifferent fixed point. The consequences of this are also very 
similar. Points which land close to the critical point tend to 
remain close to it's orbit  for a 
particularly long time. Thus even though the orbit is not a fixed or 
even a periodic point, and may even be dense, it makes sense to still 
think of nearby points as being \emph{trapped} by it for a certain 
length of time which depends on the particular rate of growth of the 
derivative.  During this time small intervals are not 
distributing themselves over the whole space as uniformly as they 
should and thus the mixing process is delayed and the rate of decay of 
correlations is correspondingly slower. 
The same argument works also in the case in which the rate of growth 
of the derivatives \( D_{n} \) is exponential. In this case the 
critical orbits behave analogously to \emph{hyperbolic repelling} 
orbits and nearby points are pushed away exponentially fast. Thus 
there is no significant loss in the rate of mixing and the decay of 
correlations is exponential in this case.

\begin{theorem}\cite{BruLuzStr03}\label{decay}
Let $f$ satisfy 
\begin{equation*}
\sum_n D_n^{-1/(2\ell-1)} < \infty \
\end{equation*}
for each critical point \( c \),
and let $\mu$
be an absolutely continuous invariant probability measure
with support $\mbox{supp}(\mu)$.
If $f$ is not renormalizable on $\mbox{supp}(\mu)$, then
$(\mbox{supp}(\mu), \mu, f)$ is mixing with the following rates:
\begin{description}
  \item [Polynomial case]
 If there exists \( C>0, \tau > 2\ell - 1 \) such that
 $$
D_n(c) \geq C n^{\tau}, \ 
 $$
for all
$c \in \mathcal C$ and $n \geq 1$,
then, for any \( \tilde \tau < \tfrac{\tau-1}{\ell-1} - 1 \), we have
\[
\mathcal C_n = \mathcal{O} (n^{-\tilde \tau})
\]
 \item [Exponential case]
 If there exist $C, \beta > 0$ such that
 $$
 D_n(c) \geq C e^{\beta n}
 $$
for all $c \in \mathcal C$ and $n\geq 1$, 
then there exist $\tilde \beta
> 0$ such
that
 \[
 \mathcal C_n = \mathcal\theta
e^{-\tilde\beta n}.
 \]
   \end{description}
\end{theorem}
 Thus the \emph{existence} of the critical point, which certainly plays a fundamental 
 role in determining many characteristics of the dynamics of these 
 maps, can be considered a bit of a red herring as far as the rate of 
 decay of correlations in concerned, i.e. it does not seem to be a 
 crucial ingredient \emph{in and of itself}.  It is more useful to keep in 
 mind the property of uniform expansion 
 outside some small neighbourhood  \( \Delta \) of the critical 
 points, and to think of \( \Delta \) in a similar way to the 
 neighbourhood of the indifferent fixed point in the previous class of 
 examples. Due to the uniform expansion outside \( \Delta \), points 
 behave in a stochastic-like away there and thus tend to fall in \( 
 \Delta \) with some definite frequency. Once they fall in \( \Delta \) 
 they remain \emph{trapped} for a certain amount of time, not in \( 
 \Delta \) itself, but in a small neighourhood of some finite number 
 of iterates of the corresponding critical orbit. If the critical 
 orbit is neutral, this slows down the mixing process, but if it is 
 exponential, this number of iterates is so small that it hardly 
 affects the mixing process. 

 We remark that exponential decay of correlations in the unimodal 
 case has been proved in \cites{You92, BalVia96} assuming exponential 
 derivative growth and a bounded recurrence 
 condition along the critical orbit, and in \cite{KelNow92} assuming 
 only the exponential growth condition. The extension of this case
 to the multimodal context, and the subexponential estimates, 
 are proved for the first time in \cite{BruLuzStr03}

\section{Viana maps}
 
Viana maps were introduced in \cite{Via97} as an example of a class of 
systems which are 
strictly \emph{not} uniformly expanding but for which the non-uniform 
expansivity condition is satisfied and, most remarkably, is
\emph{persistent} under small \( C^{3} \) perturbations, which is not 
the case for any of the examples discussed above. These maps are 
defined as skew-products on a two dimensional cylinder of the form
 \(  f: \mathbb S^1\times\mathbb R
\to \mathbb S^1\times \mathbb R \)
\[ 
f (\theta, x) = 
( \kappa \theta, x^{2}+a + \varepsilon \sin 2\pi\theta)
 \]
where \( \varepsilon \) is assumed sufficiently small and \( a \) is 
chosen so that the one-dimensional quadratic map \( x\mapsto x^{2}+a \) 
for which the critical point lands after a finite number of iterates 
onto a hyperbolic repelling periodic orbit (and thus is a \emph{good} 
parameter value and satisfies the non-uniform expansivity conditions 
as mentioned above). The map \( \kappa\theta \) is 
 taken modulo \( 2\pi \), and the constant \( \kappa \) is a 
positive integer which was required to be \( \geq 16 \) in \cite{Via97} 
although it was later shown in \cite{BuzSesTsu03} that any integer \( \geq 2 \) 
will work. The \( \sin \) function in the skew product can also be 
replaced by more general Morse functions. 

Viana proved that such class of skew-products are non-uniformly 
expanding, by showing directly 
that Lebesgue almost every point satisfies the 
non-uniform expansivity condition. It was then proved in 
\cites{Alv00,AlvVia02}
that $f$ is topologically mixing and
has a unique ergodic 
invariant measure which is absolutely continuous with respect to the 
 two-dimensional Lebesgue measure. 
As regards the rate of decay of correlation however, the situation is 
more subtle than in the classes of examples mentioned 
above. There is no particular finite set of orbits which 
to play the role of either the indifferent fixed point or the 
critical points. There is a whole curve of critical points and 
different points on this curve have different behaviour which is very 
difficult to control.  Thus the characteristics of the map which 
determine the rate of decay of correlation are less easy to pinpoint 
in a geometrical sense. Nevertheless, some geometrical structure can 
be obtained and used to obtain some estimates for the rate of decay 
of correlations.  First estimates were 
obtained in \cites{AlvLuzPin03, AlvLuzPindim1, AlvLuzPin} 
as a Corollary of a general theory of decay of 
correlations for non-uniformly expanding maps which will be discussed 
in the next section.  

\begin{theorem}\cite{AlvLuzPin}
    Viana maps have super-polynomial decay of correlations: for any \( 
    \gamma > 0 \) 
    we have
    \[ \mathcal C_{n}=  \mathcal O(n^{-\gamma}) \]
\end{theorem}

Sharper results have since been obtained in \cite{BalGou} 
by concentrating on Viana maps and taking advantage of some additional 
information which is not available in the general abstract setting, 
and also adapting some ideas from previous related 
work \cite{BalBenMau00}. 

\begin{theorem}\cite{BalGou}\label{BalGou}
    For all small enough \( \varepsilon > 0 \) there exists a 
    constant \( C_{\varepsilon} > 0 \) such that 
    Viana maps have stretched exponential decay of correlations: 
    \[ 
    \mathcal C_{n} =\mathcal O
   (e^{\sqrt n/C_{\varepsilon}})
    \]
\end{theorem}

Both results rely on some estimates from \cite{Via97} and 
the results of Theorem \ref{BalGou} 
are probably the best one can do on the basis of
those estimates.
However it is quite possible that a fresh approach to Viana maps, or 
possibly even just an improved set of estimates based on the original 
arguments of Viana, could  lead to estimates giving a faster rate of 
decay of correlations.

 \section{General theory of non-uniformly expanding maps}
\label{general theory}
 So far we have addressed the question of the rate of decay of 
 correlations in several examples of non-uniformly 
 expanding maps.  In each case, the arguments used to obtain the 
 results rely on particular features of the system.  Here we want to 
 discuss some general abstract theory of non-uniformly expanding maps. It turns 
 out that it is possible to formulate a quantitative measure of the 
 \emph{degree of non-uniformity} of the system, which 
 contains information about how close or how far the system is to 
 being uniformly expanding.  The following results show that the rate 
 of decay of correlations is closely related to this measure of the 
 non-uniformity of the expansion. They have been announced in 
 \cite{AlvLuzPin03} and are given in full details in \cite{AlvLuzPindim1} for 
 the one-dimensional case, and in \cite{AlvLuzPin} in the case of maps in 
 manifolds of arbitrary dimension. We start with the case in which \( 
 f \) is a \( C^{2} \) local diffeomorphism since the conceptual 
 picture is more straightforward. We then show how they can be 
 extended to map with critical points and even discontinuities and/or 
 points with infinite derivative. 
 
 \subsection{Non-uniformly expanding local diffeomorphisms}
 
 Let \( f: M \to M \) be a \( C^{2} \) local diffeomorphism of the 
 compact manifold \( M \) of dimension \( d\geq 1 \). We suppose that \( 
 f \) satisfies the non-uniform expansivity condition given above: 
 there exists a
    constant \( \lambda>0 \) such that 
\begin{equation}\tag{\(*\)}
    \liminf_{n\to\infty}\frac{1}{n}
    \sum_{i=0}^{n-1}\log \|Df^{-1}_{f^{i}(x)}\|^{-1}>\lambda
\end{equation}
    for almost every \( x\in M \). 
For simplicity we suppose also that \( f \) is topologically 
transitive, i.e. there exists a point \( x \) whose orbit is dense in \( 
M \). 
    Since we have no geometrical information 
    whatsoever about \( f \) we want to show that the statistical 
    properties such as the rate of decay of 
    correlations somehow depends on abstract information related to 
    the non-uniform expansivity condition only. Thus we make the 
    following 
    \begin{definition}
    For \( x \in M \), we define the \emph{expansion time function} 
    \[
\mathcal E(x) =
\min\left\{N: \frac{1}{n}
\sum_{i=0}^{n-1} \log \|Df^{-1}_{f^{i}(x)}\|^{-1} \geq \lambda/2
\ \ \forall n\geq N\right\}.
\]
\end{definition}
By condition \( (*) \) this function is defined and finite almost 
everywhere. It measures the amount of time one has to wait before the 
uniform exponential growth of the derivative kicks in. If \( \mathcal 
E(x) \) was uniformly bounded, we would essentially be in the uniformly 
expanding case. In general it will take on arbitrarily large values 
and not be defined everywhere. 
If \( \mathcal E(x) \) is \emph{large} only on a \emph{small} set of 
points, then it makes sense to think of the map as being not very 
non-uniform, whereas, if it is large on a large set of points it is 
in some sense, very non-uniform. To formalize this notion, we define 
the set 
\[ 
\Gamma_{n}=\{x\in M : \mathcal E(x) > n\}
\]
and formulate our assumptions about the \emph{degree of 
non-uniformity of \( f \)} in terms of the \emph{rate of decay of the 
measure of} \( \Gamma_{n} \).

\begin{theorem}\label{t:locdif}\cites{AlvLuzPin03, AlvLuzPin, 
AlvLuzPindim1}
    Let \( f: M\to M \) be a transitive \( C^{2} \) local diffeomorphism satisfying
    condition \( (*) \)  and suppose that there exists \( \gamma>1 \) such that
    \[
    m(\Gamma_{n}) =\mathcal O(n^{-\gamma}).
    \]
    Then there exists an absolutely continuous, \( f \)-invariant,
    probability measure \( \mu \). 
    Some finite power of \( f \) is
    mixing with respect to $\mu$ and the correlation function 
    \( \mathcal C_{n} \) for H\"older
    continuous observable on \( M \) satisfies
    \[
    \mathcal C_{n} =\mathcal O(n^{-\gamma+1 }).
    \]
 \end{theorem}
    We remark that  the existence and ergodicity of the measure \( \mu \) 
    was proved in \cite{AlvBonVia00} assuming only condition \( (*) \). 
The arguments in \cites{AlvLuzPin03, AlvLuzPindim1, AlvLuzPin} 
    give an alternative proof under the additional 
 assumption on the rate of decay of \( m(\Gamma_{n}) \).  
   We also remark that the choice of \( \lambda/2 \) in the definition 
   of the expansion time function \( \mathcal E(x) \) is fairly 
   arbitrary and does not  affect the asymptotic rate estimates. Any 
   positive number smaller than \( \lambda \) would yield the same 
   results. 
   
   \subsection{Critical points and discontinuities}
   The arguments used in the proof of Theorem \ref{t:locdif} apply to 
   maps which fail to be local diffeomorphism on some zero measure set 
   \( \mathcal C \) satisfying some simple non-degeneracy conditions. 
   Points in \( \mathcal C \) maybe the the higher dimensional 
   analogue of critical points (i.e. points at which derivative is 
   degenerate), or may be points of discontinuities for \( f \), or 
   points at which the map is not differentiable and for which the 
   derivative blows up to infinity. Remarkably, all these cases are 
   dealt with, as ``problematic'' points, in the same way and 
 need to satisfy the same conditions which are just the natural 
   generalization of the non-degeneracy (or non-flatness) condition 
   for critical points of one-dimensional maps. 
   
   \begin{definition}
    The \emph{critical set} \( \mathcal C \subset M \) is
    \emph{non-degenerate} if \( m(\mathcal C)=0 \) 
    and there is a constant $\beta>0$
    such that for every $x\in M\setminus\mathcal C$ we have
    $\dist(x,\mathcal C)^{\beta}\lesssim \|Df_{x}v\|/\|v\|\lesssim
    \dist(x,\mathcal C)^{-\beta}$ for all $v\in T_x M$, and the functions \(
    \log\det Df \) and \( \log \|Df^{-1}\| \) are \emph{locally Lipschitz}
    with Lipschitz constant \( \lesssim \dist(x,\mathcal C)^{-\beta}\).
\end{definition}

These are geometrical conditions which have nothing to do with the 
dynamics. We also need to assume some dynamical conditions concerning 
the rate of recurrence of typical points near the critical set. 
We let \( d_{\delta}(x,\mathcal C) \) denote the \( \delta \)-\emph{truncated}
distance from \( x \) to \( \mathcal C \) defined as \(
d_{\delta}(x,\mathcal C) = d(x,\mathcal C) \) if \( d(x,\mathcal C)
\leq \delta\) and \( d_{\delta}(x,\mathcal C) =1 \) otherwise.
\begin{definition}
We say that \( f \) satisfies the property of 
{\em subexponential
recurrence} to the critical set if
for any $\epsilon>0$ there exists $\delta>0$ such that for Lebesgue
almost every $x\in M$
\begin{equation} \label{e.faraway1}\tag{\( ** \)}
    \limsup_{n\to+\infty}
\frac{1}{n} \sum_{j=0}^{n-1}-\log \dist_\delta(f^j(x),\mathcal C)
\le\epsilon.
\end{equation}
\end{definition}
Again, we want to differentiate between 
different degrees of recurrence in a similar way to the way we 
differentiated between different degrees of non-uniformity of the 
expansion. 
\begin{definition}
       For \( x \in M \), we define the \emph{recurrence time function} 
\[
\textstyle{
\mathcal R(x) = \min\left\{N\ge 1: \frac{1}{n} \sum_{j=0}^{n-1} -\log
\dist_\delta(f^j(x),\mathcal C) \leq 2\varepsilon, \forall n\geq N\right\}
}
\]
\end{definition}
Then, for a map satisfying both conditions \( (*) \) and \( (**) \) we 
let 
\[ 
\Gamma_{n} = \{x: \mathcal E(x) > n \text{ or
} \mathcal R(x) > n\}
\]
\begin{theorem}\cites{AlvLuzPin03, AlvLuzPin, 
AlvLuzPindim1}
      Let \( f:M\to M \) be a transitive \( C^{2} \) local
 diffeomorphism outside a non-degenerate critical set \( \mathcal C 
 \),
 satisfying conditions \( (*) \) and \( (**) \). Suppose that there exists \(
 \gamma>0 \) such that 
 \[
 |\Gamma_n| = \mathcal O(n^{-\gamma}).  
 \] 
     Then there exists an absolutely continuous, \( f \)-invariant,
    probability measure \( \mu \).
    Some finite power of \( f \) is
    mixing with respect to $\mu$ and 
    for any H\"older continuous function \( \varphi, \psi \) on \(M \)
we have
    \[
    \mathcal C_{n} = \mathcal O(n^{-\gamma+1 }).
    \]
    \end{theorem}
We remark that although condition \( (**) \) might appear to be a very technical 
condition, it is actually quite natural and in fact \emph{almost} necessary. 
Indeed, suppose that an absolutely continuous invariant 
measure \( \mu \) did exist for \( f \). Then, a simple application of 
Birkhoff's Ergodic theorem implies that condition \( (**) \) is 
equivalent to the integrability condition 
\[ 
\int_{M} |\log \dist_\delta(x,\mathcal S)| d\mu < \infty
\]
which is simply saying that the invariant measure 
does not give too 
much weight to a neighbourhood of the discontinuity set.

\section{Concluding remarks}
    
\subsection{What causes slow decay of correlations ?}

The results on maps with indifferent fixed points and those on one 
dimensional maps with critical points suggest that slow decay of 
correlation is caused, literally, by a slowing down of the dynamics due 
to some indifferent fixed or periodic point, or even some indifferent 
non-periodic orbit. In these cases, the responsible orbits can be 
identified exactly and the mechanism through which they slow down the 
mixing process can be described quite explicitly and even quantified 
in terms of the degree of ``neutrality'' of the orbit. In higher 
dimensional cases where it is impossible to identify specific orbits 
with neutral behaviour it is much more difficult to obtain optimal 
estimates. It is not known if Viana maps have any neutral orbits and 
the current estimates for decay of correlations are based on estimates 
which may not be optimal. The results of \cites{AlvLuzPin03, AlvLuzPin, 
AlvLuzPindim1} described 
in the last section develop a connection between the rate of decay of 
correlations and the  rate at which orbits start to exhibit 
exponential growth of the derivative. However this theory still does 
resolve the issue of what geometrical or other characteristics of 
the system will cause this rate to be of a certain kind rather than 
another. In view of the above discussion, it seems reasonable to 
conjecture in the first instance that the difference between 
exponential and subexponential behaviour lies in the absence or 
presence of a neutral orbit of some kind. 

To formalize this notion we recall a few standard notions of Ergodic 
Theory. We suppose that \( f:M \to M \) is a \( C^{2} \) map of the 
compact manifold \( M \) of dimension \( d\geq 1 \), and we let \( 
\mathcal M_{inv} \) denote the space of all probability invariant 
measures \( \mu \) on \( M \) which satisfy the integrability 
condition 
\[ 
\int_{M} \log \|Df_{x}\| d\mu < \infty.
\]
Then we can apply a version of Oseledet's Theorem for non-invertible maps 
which says 
that there exist constants \( \lambda_{1}, \ldots, \lambda_{k} \) 
with \( k\leq d \), and 
a measurable decomposition \( T_{x}M = E^{1}_{x}\oplus\dots\oplus E^{k}_{x} \)
of the tangent bundle over \( M \) such that the decomposition is 
invariant by the derivative and such that for all \( j=1,\ldots, k \) 
and for all non zero vectors \( v^{(j)}\in 
E^{j}_{x} \) we have 
\[ 
\lim_{n\to\infty}\frac{1}{n}\sum_{i=0}^{n-1}\log \|Df^{n}_{x}(v^{(j)})\| = 
\lambda_{j}.
\]
The constants \( \lambda_{1}, \ldots, \lambda_{k} \) are called the 
\emph{Lyapunov exponents} associated to the measure \( \mu \). The 
set of Lyapunov exponents associated to \( \mu \) is sometimes called 
the \emph{Lyapunov spectrum} of \( \mu \) and denoted by \emph{Lyap}(\( 
\mu \)). The standard terminology is to call \( \mu \) hyperbolic if 
it has no associated zero Lyapunov exponents. 

\begin{definition}\label{expanding measure}
    \( \mu \in \mathcal M_{inv}\) is
\emph{expanding} if  every Lyapunov exponent associated to \( \mu \) 
is \( > 0 \). 
\end{definition}

Notice that if \( f \) satisfies the non-uniform expansivity condition 
\( (*) \) defined above, then the corresponding absolutely continuous 
invariant measure \( \mu \) (which exists by \cite{AlvBonVia00}) is 
expanding in the sense of definition \ref{expanding measure}. 
Conversely, if \( \mu \) is an expanding absolutely continuous 
invariant probability measure, then it clearly satisfies condition 
\( (*) \). Thus these notions characterize 
the class of (non-uniformly) expanding maps. Notice however that 
in general there are infinitely many invariant measures, for example 
the Dirac measures on fixed and periodic points. Even if the 
absolutely continuous invariant measure is expanding, there may be 
other singular measures which are not, i.e. have a zero Lyapunov 
exponent. To distinguish between these two cases we formulate the 
following

\begin{definition}
    We say that \( f \) is \emph{totally expanding} if all measures 
    in \( \mathcal M_{inv} \) are \emph{uniformly} expanding in the sense that  \( \exists \ 
    \lambda > 0 \) such that every Lyapunov exponent of every measure \( 
    \mu\in\mathcal M_{inv} \) is \( \geq \lambda \). 
\end{definition}
 
A totally expanding map is in some sense the strongest possible 
version of a non-uniformly expanding map without necessarily being 
uniformly expanding: a kind of ``uniformly'' 
non-uniformly expanding !
 Uniformly expanding maps are clearly totally expanding.  
 Remarkably, if \( f \) is a \( C^{1} \) local diffeomorphisms then 
 the two notions are equivalent: completely 
 expanding implies uniformly expanding 
 \cites{AlvAraSau, Cao}.  
 However general non-uniformly expanding maps are certainly not 
 necessarily totally expanding. 
The examples of maps with indifferent fixed points of Section 
 \ref{almost expanding} are not totally expanding, 
 as the Dirac measure on the 
 indifferent fixed point is not expanding. In the context of 
 one-dimensional smooth maps with critical points it is known that in 
 the unimodal case exponential growth of the derivative along the 
 critical orbit (the Collet-Eckmann condition) 
 implies uniform hyperbolicity on periodic orbits 
 \cite{Now88} which in turn implies total expansivity 
 \cite{BruKel98}. Conversely total expansivity clearly implies uniform 
 hyperbolicity on periodic orbits which implies exponential growth of 
 the derivative along the critical orbit  \cite{NowSan98}.
 Moreover \cite{NowSan98} closes a chain of implications which 
 finally imply that in the unimodal case, total hyperbolicity is 
 equivalent to exponential decay of correlations. This suggests the 
 following
 
 \begin{conjecture*}
     Suppose \( f: M \to M \) is non-uniformly expanding with 
     absolutely continuous invariant measure \( \mu \). 
     Then \( f \) 
     has exponential decay of correlations if and only if it is 
     totally expanding
  \end{conjecture*}
  
  The idea is that any invariant measure with a zero Lyapunov 
  exponent causes a slowing down effect analogous to that caused by 
  an indifferent   fixed point or an indifferent critical orbit. This 
  conjecture is supported by known examples as described above, 
  but at the moment 
  it is not clear how it could be tackled in general. Even the implication in 
  just one of the directions would interesting
  and  support the intuitive picture of the cause of 
  slow decay of correlation given in the discussion. 
 
We remark that the assumption of non-uniform expansion is 
crucial here. There are several examples of systems which 
have exponential decay of correlations but clearly have invariant 
measures with zero Lyapunov exponents, e.g. 
partially hyperbolic maps or maps obtained as time-1 maps of certain 
flows \cites{Dol98a, Dol98b, Dol00}. These examples however are 
not non-uniformly expanding, and are generally \emph{partially hyperbolic} 
which means that there are two continuous subbundles such that the 
derivative restricted to one subbundle has very good expanding 
properties or contracting properties and the other subbundle has the 
zero Lyapunov exponents. For reasons which are not at all clear, this 
might be \emph{better} from the point of view of decay of correlations 
than a situation in which all the Lyapunov exponents of the absolutely 
continuous measure are positive but there is some \emph{embedded} 
singular measure with zero Lyapunov exponent slowing down the mixing 
process. Certainly there is still a lot to be understood on this topic.


\subsection{Arguments and techniques}
The existence of absolutely continuous invariant measures can, in 
many cases, be proved by fairly direct geometric arguments. Estimates 
on the rate of decay of correlations, on the other hand usually
require considerably more sophisticated arguments. In the case of 
expanding maps, these generally involve an abstract 
functional-analytic or probabilistic framework in which particular 
geometric characteristics of the system (expansion, smoothness) 
are used as ingredients.  

The pioneering work of Sinai, Ruelle and Bowen on uniformly 
hyperbolic systems (with uniformly expanding systems essentially as 
a special case), introduced 
the basic idea of approaching the problem via the Perron-Frobenius 
operator. This is an operator on a suitable function space of 
appropriate \emph{densities} with the property that any fixed point 
for the operator is the density of an invariant absolutely continuous 
measure. Thus, various functional-analytic techniques can be brought 
to bear on the problem of the existence of such a fixed point and on 
the speed of convergence of arbitrary densities to the fixed points 
which turns out to be closely related to the rates of decay of 
correlations. In particular, a \emph{spectral gap} in the spectrum of 
the operator implies exponential decay of correlations. Despite the 
great success of such an approach in dealing with various classes of 
systems, it has the intrinsic limitation of producing results which 
are necessarily exponential: if there is no spectral gap then one 
cannot deduce any other rate for the decay. 
Another 
functional-analytic method has been introduced in \cite{Liv95} to 
deal with  maps with discontinuities. This 
still involves a direct study of the Perron-Frobenius operator but 
using the so-called \emph{Birkhoff metric} and the notion of 
\emph{invariant cones}. This method appears to have more in-built 
flexibility and in certain cases allows better estimates for the 
actual constants and exponents involved in decay of correlations. 
Moreover it has been adapted in \cite{Mau01b} to deal with systems with 
subexponential decay. 

More recently, a quite different approach has been introduced in 
\cites{You98, BreFerGal99, You99} which relies on a probabilistic \emph{coupling} 
argument. As in the functional-analytic approach, it depends 
ultimately on the fact that the rate of decay  correlations is 
related to the rate at which  arbitrary absolutely continuous 
measures with densities satisfying some regularity conditions converge 
to the invariant measure under the dynamics. However, in this case 
the conclusions do not ultimately depend on a spectral estimate, but on
more direct geometric estimates. 
This approach has proved to be very flexible and far reaching in its 
scope and underlies several recent results such as those in 
\cites{BruLuzStr03, AlvLuzPindim1, AlvLuzPin03} as well as those in 
the pioneering papers \cites{You98, You99}.  It has also 
proved successful in generalizing known results to observables which 
satisfy significantly weaker summability conditions on the modulus of 
continuity of the derivative, rather than the usual H\"older 
continuity conditions \cite{Lyn03}. It is beyond the scope of this 
note to enter into a more detailed discussion of the various techniques. 
See the excellent and comprehensives texts \cites{Bal00, Via} for 
detailed discussions of the functional-analytic methods in particular, 
 and \cite{You99} for the coupling method. 

  \subsection{Lower bounds for decay of correlations}
Finally we make the important observation that all the estimates given 
above have been upper bounds for the decay of correlations. 
An important question for a full understanding of the phenomenon of 
decay of correlations is that of whether these are actually lower 
bounds as well. It turns out that lower bounds require some 
sophisticated arguments which are not always directly related to the 
arguments used to prove upper bounds. There are not many results in 
this direction although there have been some important recent 
developments in \cites{Sar02, Gou}. 
 
 \subsection{Hyperbolicity}
 There is also extensive research work on higher dimensional systems 
 which exhibit some degree of \emph{hyperbolicity}, 
 a combination of expansion and contraction in 
 different directions. As well as the natural analogues of the 
 uniformly and non-uniformly expanding systems, uniformly and 
 non-uniformly hyperbolic systems, there are other categories such as 
 partially hyperbolic or projectively hyperbolic. These definitions 
 all try to distinguish different possible kinds of hyperbolicity 
 both in terms of how the expansion and contraction estimates behave 
 within certain contracting and expanding subbundles of the tangent 
 space, as well as how these subbundles are related to each other. 
 Again we refer to \cites{Bal00, Via} for more details and additional 
 references.


\begin{bibsection}[Bibliography]
\begin{biblist}
    \bib{Aar97}{book}{
  author={Aaronson, Jon},
  title={An introduction to infinite ergodic theory},
  volume={50},
  series={Mathematical surveys and monographs},
  publisher={AMS},
}
\bib{Adl73}{article}{
  author={Adler, Roy L.},
  title={$F$-expansions revisited},
  booktitle={Recent advances in topological dynamics (Proc. Conf., Yale Univ., 
  New Haven, Conn., 1972; in honor of Gustav Arnold Nedlund)},
  pages={1\ndash 5. Lecture Notes in Math., Vol. 318},
  publisher={Springer},
  place={Berlin},
  date={1973},
}
\bib{Alv00}{article}{
  author={Alves, Jos{\'e} Ferreira},
  title={SRB measures for non-hyperbolic systems with multidimensional expansion},
  language={English, with English and French summaries},
  journal={Ann. Sci. \'Ecole Norm. Sup. (4)},
  volume={33},
  date={2000},
  number={1},
  pages={1\ndash 32},
}
\bib{AlvAraSau}{article}{
  author={Alves, Jos{\'e} Ferreira},
  author={Araujo, Vitor},
  author={Saussol, Benoit},
  title={On the uniform hyperbolicity of some nonuniformly hyperbolic systems},
  status={preprint},
  year={2002},
}
\bib{AlvBonVia00}{article}{
  author={Alves, Jos{\'e} F.},
  author={Bonatti, Christian},
  author={Viana, Marcelo},
  title={SRB measures for partially hyperbolic systems whose central 
  direction is mostly expanding},
  journal={Invent. Math.},
  volume={140},
  date={2000},
  number={2},
  pages={351\ndash 398},
}
\bib{AlvLuzPin03}{article}{
  author={Alves, Jos{\'e} F.},
  author={Luzzatto, Stefano},
  author={Pinheiro, Vilton},
  title={Markov structures and decay of correlations for non-uniformly 
  expanding maps on compact manifolds of arbitrary dimension.},
  journal={Electronic Research Announcements of the AMS},
  year={2003},
  eprint={\href{http://www.ams.org/journal-getitem?pii=S1079-6762-03-00106-9}
  {http://www.ams.org/journal-getitem?pii=S1079-6762-03-00106-9}},
}
\bib{AlvLuzPindim1}{article}{
  author={Alves, Jos{\'e} F.},
  author={Luzzatto, Stefano},
  author={Pinheiro, Vilton},
  title={Lyapunov exponents and rates of mixing for one-dimensional maps.},
  status={to appear},
  journal={Erg. Th. \& Dyn. Syst.},
}
\bib{AlvLuzPin}{article}{
  author={Alves, Jos{\'e} F.},
  author={Luzzatto, Stefano},
  author={Pinheiro, Vilton},
  title={Markov structures and decay of correlations for non-uniformly expanding maps.},
  status={preprint},
  eprint={\href{http://front.math.ucdavis.edu/math.DS/0205191}
  {http://front.math.ucdavis.edu/math.DS/0205191}},
}
\bib{AlvVia02}{article}{
  author={Alves, Jos{\'e} F.},
  author={Viana, Marcelo},
  title={Statistical stability for robust classes of maps with non-uniform expansion},
  journal={Ergodic Theory Dynam. Systems},
  volume={22},
  date={2002},
  number={1},
  pages={1\ndash 32},
}
\bib{Ave68}{article}{
  author={Avez, Andr{\'e}},
  title={Propri\'et\'es ergodiques des endomorphisms dilatants des vari\'et\'es compactes},
  language={French},
  journal={C. R. Acad. Sci. Paris S\'er. A-B},
  volume={266},
  date={1968},
  pages={A610\ndash A612},
}
\bib{Bak97}{book}{
  author={Bak, Per},
  title={How nature works},
  publisher={OUP},
  year={1997},
}
\bib{Bal00}{book}{
  author={Baladi, Viviane},
  title={Positive transfer operators and decay of correlations},
  series={Advanced Series in Nonlinear Dynamics},
  volume={16},
  publisher={World Scientific Publishing Co. Inc.},
  place={River Edge, NJ},
  date={2000},
  pages={x+314},
}
\bib{BalBenMau00}{article}{
  author={Baladi, Viviane},
  author={Benedicks, Michael},
  author={Maume-Deschamps, V{\'e}ronique},
  title={Decay of random correlation functions for unimodal maps},
  note={XXXI Symposium on Mathematical Physics (Toru\'n, 1999)},
  journal={Rep. Math. Phys.},
  volume={46},
  date={2000},
  number={1-2},
  pages={15\ndash 26},
}
\bib{BalGou}{article}{
  author={Baladi, Viviane},
  author={Gou\"ezel, S\'ebastien},
  title={Stretched exponential bounds for the correlations of the Alves-Viana skew product},
  status={preprint},
}
\bib{BalVia96}{article}{
  author={Baladi, Viviane},
  author={Viana, Marcelo},
  title={Strong stochastic stability and rate of mixing for unimodal maps},
  journal={Ann. Sci. \'Ecole Norm. Sup. (4)},
  volume={29},
  date={1996},
  number={4},
  pages={483\ndash 517},
}
\bib{BenCar85}{article}{
  author={Benedicks, M.},
  author={Carleson, L.},
  title={On iterations of $1-ax^2$ on $(-1,1)$},
  date={1985},
  journal={Ann. of Math.},
  volume={122},
  pages={1\ndash 25},
}
\bib{Bow75}{book}{
  author={Bowen, Rufus},
  title={Equilibrium states and the ergodic theory of Anosov diffeomorphisms},
  note={Lecture Notes in Mathematics, Vol. 470},
  publisher={Springer-Verlag},
  place={Berlin},
  date={1975},
  pages={i+108},
}
\bib{Bow77}{article}{
  author={Bowen, Rufus},
  title={Bernoulli maps of the interval},
  journal={Israel J. Math.},
  volume={28},
  date={1977},
  number={1-2},
  pages={161\ndash 168},
}
\bib{Bow79}{article}{
  author={Bowen, Rufus},
  title={Invariant measures for Markov maps of the interval},
  note={With an afterword by Roy L. Adler and additional comments by Caroline Series},
  journal={Comm. Math. Phys.},
  volume={69},
  date={1979},
  number={1},
  pages={1\ndash 17},
}
\bib{Bre99}{article}{
  author={Bressaud, Xavier},
  title={Subshifts on an infinite alphabet},
  journal={Ergodic Theory Dynam. Systems},
  volume={19},
  date={1999},
  number={5},
  pages={1175\ndash 1200},

}
\bib{BreFerGal99}{article}{ 
author={Bressaud, Xavier},
author={Fern{\'a}ndez, Roberto}, 
author={Galves, Antonio}, 
title={Decay of
correlations for non-H\"olderian dynamics.  A coupling approach},
journal={Electron.  J. Probab.}, 
volume={4}, date={1999}, 
pages={no.  3, 19
pp.  (electronic)}, 
}

\bib{BruKel98}{article}{
  author={Bruin, Henk},
  author={Keller, Gerhard},
  title={Equilibrium states for $S$-unimodal maps},
  journal={Ergodic Theory Dynam. Systems},
  volume={18},
  date={1998},
  number={4},
  pages={765\ndash 789},
}
\bib{BruLuzStr03}{article}{
  author={Bruin, Henk},
  author={Luzzatto, Stefano},
  author={van Strien, Sebastian},
  title={Decay of correlations in one-dimensional dynamics},
  journal={Ann. Ec. Norm. Sup.},
  year={2003},
}
\bib{BruStr01}{article}{
  author={Bruin, Henk},
  author={van Strien, Sebastian},
  title={Existence of absolutely continuous invariant probability measures for multimodal maps},
  booktitle={Global analysis of dynamical systems},
  pages={433\ndash 447},
  publisher={Inst. Phys.},
  place={Bristol},
  date={2001},
}
\bib{BruStrShe03}{article}{
  author={Bruin, Henk},
  author={van Strien, Sebastian},
  author={Shen, Weixiao},
  title={Invariant measures exist without a growth conditions},
  status={preprint},
}
\bib{Buz00a}{article}{
  author={Buzzi, J{\'e}r{\^o}me},
  title={Absolutely continuous invariant probability measures for arbitrary expanding piecewise $\mathbf R$-analytic mappings of the plane},
  journal={Ergodic Theory Dynam. Systems},
  volume={20},
  date={2000},
  number={3},
  pages={697\ndash 708},
}
\bib{Buz01a}{article}{
  author={Buzzi, J{\'e}r{\^o}me},
  title={No or infinitely many a.c.i.p.\ for piecewise expanding $C\sp r$ maps in higher dimensions},
  journal={Comm. Math. Phys.},
  volume={222},
  date={2001},
  number={3},
  pages={495\ndash 501},
}
\bib{Buz01c}{article}{
  author={Buzzi, J{\'e}r{\^o}me},
  title={Thermodynamical formalism for piecewise invertible maps: absolutely continuous invariant measures as equilibrium states},
  booktitle={Smooth ergodic theory and its applications (Seattle, WA, 1999)},
  language={English, with English and French summaries},
  series={Proc. Sympos. Pure Math.},
  volume={69},
  pages={749\ndash 783},
  publisher={Amer. Math. Soc.},
  place={Providence, RI},
  date={2001},
}
\bib{Buz99a}{article}{
  author={Buzzi, J.},
  title={Absolutely continuous invariant measures for generic multi-dimensional piecewise affine expanding maps},
  note={Discrete dynamical systems},
  journal={Internat. J. Bifur. Chaos Appl. Sci. Engrg.},
  volume={9},
  date={1999},
  number={9},
  pages={1743\ndash 1750},
}
\bib{BuzKel01}{article}{
  author={Buzzi, J{\'e}r{\^o}me},
  author={Keller, Gerhard},
  title={Zeta functions and transfer operators for multidimensional piecewise affine and expanding maps},
  journal={Ergodic Theory Dynam. Systems},
  volume={21},
  date={2001},
  number={3},
  pages={689\ndash 716},
}
\bib{BuzMau03}{article}{
  author={Buzzi, J{\'e}r{\^o}me},
  author={Maume, Veronique},
  title={Decay of correlations for piecewise invertible maps in higher dimensions},
  status={Preprint},
}
\bib{BuzPacSch01}{article}{
  author={Buzzi, J{\'e}r{\^o}me},
  author={Paccaut, Fr{\'e}d{\'e}ric},
  author={Schmitt, Bernard},
  title={Conformal measures for multidimensional piecewise invertible maps},
  journal={Ergodic Theory Dynam. Systems},
  volume={21},
  date={2001},
  number={4},
  pages={1035\ndash 1049},
}
\bib{BuzSar}{article}{
  author={J{\'e}r{\^o}me Buzzi},
  author={Omri Sarig},
  title={Uniqueness of equilibrium measures for countable Markov shifts and 
  multi-dimensional piecewise expanding maps},
  status={Preprint},
}
\bib{BuzSesTsu03}{article}{
  author={Buzzi, Jerome},
  author={Sester, 0},
  author={Tsujii, Masato},
  title={ Weakly expanding skew-products of quadratic maps},
  journal={Ergod. Th. adn Dyn. Syst.},
  year={2003},
}
\bib{Cao}{article}{
  author={Cao, Yongluo},
  title={Positive Lyapunov exponents and uniform hyperbolicity},
  status={preprint},
  year={2003},
}
\bib{ColEck83}{article}{
  author={Collet, P.},
  author={Eckmann, J.-P.},
  title={Positive Liapunov exponents and absolute continuity for maps of the interval},
  journal={Ergodic Theory Dynam. Systems},
  volume={3},
  date={1983},
  number={1},
  pages={13\ndash 46},
}
\bib{Cow02}{article}{
  author={Cowieson, William J.},
  title={Absolutely continuous invariant measures for most piecewise smooth expanding maps},
  journal={Ergodic Theory Dynam. Systems},
  volume={22},
  date={2002},
  number={4},
  pages={1061\ndash 1078},
}
\bib{Dol00}{article}{
  author={Dolgopyat, Dmitry},
  title={On dynamics of mostly contracting diffeomorphisms},
  journal={Comm. Math. Phys.},
  volume={213},
  date={2000},
  number={1},
  pages={181\ndash 201},
}
\bib{Dol98a}{article}{
  author={Dolgopyat, Dmitry},
  title={On decay of correlations in Anosov flows},
  journal={Ann. of Math. (2)},
  volume={147},
  date={1998},
  number={2},
  pages={357\ndash 390},
}
\bib{Dol98b}{article}{
  author={Dolgopyat, Dmitry},
  title={Prevalence of rapid mixing in hyperbolic flows},
  journal={Ergodic Theory Dynam. Systems},
  volume={18},
  date={1998},
  number={5},
  pages={1097\ndash 1114},
}
\bib{Gel59}{article}{
  author={Gel{\cprime }fond, A. O.},
  title={A common property of number systems},
  language={Russian},
  journal={Izv. Akad. Nauk SSSR. Ser. Mat.},
  volume={23},
  date={1959},
  pages={809\ndash 814},
}
\bib{GorBoy89}{article}{
  author={G{\'o}ra, P.},
  author={Boyarsky, A.},
  title={Absolutely continuous invariant measures for piecewise expanding 
  $C\sp 2$ transformation in ${\bf R}\sp N$},
  journal={Israel J. Math.},
  volume={67},
  date={1989},
  number={3},
  pages={272\ndash 286},
}
\bib{GorSch89}{article}{
  author={G{\'o}ra, P.},
  author={Schmitt, B.},
  title={Un exemple de transformation dilatante et $C\sp 1$ par morceaux de l'intervalle, sans probabilit\'e absolument continue invariante},
  language={French, with English summary},
  journal={Ergodic Theory Dynam. Systems},
  volume={9},
  date={1989},
  number={1},
  pages={101\ndash 113},
}
\bib{Gou}{article}{
  author={Gou\"ezel, S\'ebastien},
  title={Sharp polynomial estimates for the decay of correlations},
  status={preprint},
}
\bib{GraSwi97}{article}{
  author={Graczyk, Jacek},
  author={{\'S}wiatek, Grzegorz},
  title={Generic hyperbolicity in the logistic family},
  journal={Ann. of Math. (2)},
  volume={146},
  date={1997},
  number={1},
  pages={1\ndash 52},
}
\bib{Hal44}{article}{
  author={Halmos, Paul R.},
  title={In general a measure preserving transformation is mixing},
  journal={Ann. of Math. (2)},
  volume={45},
  date={1944},
  pages={786\ndash 792},
}
\bib{HofKel82}{article}{
  author={Hofbauer, Franz},
  author={Keller, Gerhard},
  title={Ergodic properties of invariant measures for piecewise monotonic transformations},
  journal={Math. Z.},
  volume={180},
  date={1982},
  number={1},
  pages={119\ndash 140},
}
\bib{Hol03}{article}{
  author={Holland, Mark},
  title={Slowly mixing systems and intermittency maps},
  status={preprint},
  year={2003},
}
\bib{Hol57}{article}{
  author={Holladay, John C.},
  title={On the existence of a mixing measure},
  journal={Proc. Amer. Math. Soc.},
  volume={8},
  date={1957},
  pages={887\ndash 893},
}
\bib{Hu01}{article}{
  author={Hu, Huyi},
  title={Statistical properties of some almost hyperbolic systems},
  booktitle={Smooth ergodic theory and its applications (Seattle, WA, 1999)},
  series={Proc. Sympos. Pure Math.},
  volume={69},
  pages={367\ndash 384},
  publisher={Amer. Math. Soc.},
  place={Providence, RI},
}
\bib{HuYou95}{article}{
  author={Hu, Hu Yi},
  author={Young, Lai-Sang},
  title={Nonexistence of SBR measures for some diffeomorphisms that are ``almost Anosov''},
  journal={Ergodic Theory Dynam. Systems},
  volume={15},
  date={1995},
  number={1},
  pages={67\ndash 76},
}
\bib{Iso99}{article}{
  author={Isola, Stefano},
  title={Renewal sequences and intermittency},
  journal={J. Statist. Phys.},
  volume={97},
  date={1999},
  number={1-2},
  pages={263\ndash 280},
}
\bib{Jak81}{article}{
  author={Jakobson, M.~V.},
  title={Absolutely continuous invariant measures for one\ndash parameter families of one\ndash dimensional maps},
  date={1981},
  journal={Comm. Math. Phys.},
  volume={81},
  pages={39\ndash 88},
}
\bib{Jen98}{book}{
  author={Jensen, Henrik Jeldtoft},
  title={Self-organized criticality : emergent complex behavior in physical and biological systems},
  publisher={CUP},
  year={1998},
}
\bib{KacKes58}{article}{
  author={Kac, M.},
  author={Kesten, Harry},
  title={On rapidly mixing transformations and an application to continued fractions},
  journal={Bull. Amer. Math. Soc. 64 (1958), 283--283; correction},
  volume={65},
  date={1958},
  pages={67},
}
\bib{Kel79}{article}{
  author={Keller, Gerhard},
  title={Ergodicit\'e et mesures invariantes pour les transformations dilatantes par morceaux d'une r\'egion born\'ee du plan},
  language={French, with English summary},
  journal={C. R. Acad. Sci. Paris S\'er. A-B},
  volume={289},
  date={1979},
  number={12},
  pages={A625\ndash A627},
}
\bib{Kel80}{article}{
  author={Keller, Gerhard},
  title={Un th\'eor\`eme de la limite centrale pour une classe de transformations monotones par morceaux},
  language={French, with English summary},
  journal={C. R. Acad. Sci. Paris S\'er. A-B},
  volume={291},
  date={1980},
  number={2},
  pages={A155\ndash A158},
}
\bib{KelNow92}{article}{
  author={Keller, Gerhard},
  author={Nowicki, Tomasz},
  title={Spectral theory, zeta functions and the distribution of periodic 
  points for Collet-Eckmann maps},
  journal={Comm. Math. Phys.},
  volume={149},
  date={1992},
  number={1},
  pages={31\ndash 69},
}
\bib{Koz03}{article}{
  author={KOZLOVSKI, O.S.},
  title={Axiom A maps are dense in the space of unimodal maps in the \( C^k \) topology},
  journal={Annals of Math.},
  year={2003},
  volume={157},
  number={1},
  pages={1\ndash 43},
}
\bib{KrzSzl69}{article}{
  author={Krzy{\.z}ewski, K.},
  author={Szlenk, W.},
  title={On invariant measures for expanding differentiable mappings},
  journal={Studia Math.},
  volume={33},
  date={1969},
  pages={83\ndash 92},
}
\bib{LSV98}{article}{
  author={Liverani, Carlangelo},
  author={Saussol, Benoit},
  author={Vaienti, Sandro},
  title={Conformal measure and decay of correlation for covering weighted systems},
  journal={Ergodic Theory Dynam. Systems},
  volume={18},
  date={1998},
  number={6},
  pages={1399\ndash 1420},
}
\bib{Las73}{article}{
  author={Lasota, A.},
  title={Invariant measures and functional equations},
  journal={Aequationes Math.},
  volume={9},
  date={1973},
  pages={193\ndash 200},
}
\bib{LasYor73}{article}{
  author={Lasota, A.},
  author={Yorke, James A.},
  title={On the existence of invariant measures for piecewise monotonic transformations},
  journal={Trans. Amer. Math. Soc.},
  volume={186},
  date={1973},
  pages={481\ndash 488 (1974)},
}
\bib{Liv95a}{article}{
  author={Liverani, Carlangelo},
  title={Decay of correlations for piecewise expanding maps},
  journal={J. Statist. Phys.},
  volume={78},
  date={1995},
  number={3-4},
  pages={1111\ndash 1129},
}
\bib{Liv95}{article}{
  author={Liverani, Carlangelo},
  title={Decay of correlations},
  journal={Ann. of Math. (2)},
  volume={142},
  date={1995},
  number={2},
  pages={239\ndash 301},
}
\bib{Luz00}{article}{
  author={Luzzatto, Stefano},
  title={Bounded recurrence of critical points and Jakobson's theorem},
  booktitle={The Mandelbrot set, theme and variations},
  series={London Math. Soc. Lecture Note Ser.},
  volume={274},
  pages={173\ndash 210},
  publisher={Cambridge Univ. Press},
  place={Cambridge},
  date={2000},
}
\bib{LuzTuc99}{article}{
  author={Luzzatto, Stefano},
  author={Tucker, Warwick},
  title={Non-uniformly expanding dynamics in maps with singularities and criticalities},
  journal={Inst. Hautes \'Etudes Sci. Publ. Math.},
  number={89},
  date={1999},
  pages={179\ndash 226 (1999)},
}
\bib{LuzVia00}{article}{
  author={Luzzatto, Stefano},
  author={Viana, Marcelo},
  title={Positive Lyapunov exponents for Lorenz-like families with criticalities},
  language={English, with English and French summaries},
  note={G\'eom\'etrie complexe et syst\`emes dynamiques (Orsay, 1995)},
  journal={Ast\'erisque},
  number={261},
  date={2000},
  pages={xiii, 201\ndash 237},
}
\bib{Lyn03}{article}{
  author={Lynch, Vincent},
  title={ Decay of correlations for 
  non-H\"older continuous observables},
  status={preprint},
}
\bib{Lyu97}{article}{
  author={Lyubich, Mikhail},
  title={Dynamics of quadratic polynomials. I, II},
  journal={Acta Math.},
  volume={178},
  date={1997},
  number={2},
  pages={185\ndash 247, 247\ndash 297},
}
\bib{Man85}{article}{
  author={Ma{\~n}{\'e}, Ricardo},
  title={Hyperbolicity, sinks and measure in one-dimensional dynamics},
  journal={Comm. Math. Phys.},
  volume={100},
  date={1985},
  number={4},
  pages={495\ndash 524},
}
\bib{Man87}{book}{
  author={Ricardo Ma\~{n}\'{e}},
  title={Ergodic theory and differentiable dynamics},
  publisher={Springer-Verlag},
  year={1987},
}
\bib{ManPom80}{article}{
  author={Manneville, Paul},
  author={Pomeau, Yves},
  title={Intermittent transition to turbulence in dissipative dynamical systems},
  journal={Comm. Math. Phys.},
  volume={74},
  date={1980},
  number={2},
  pages={189\ndash 197},
}
\bib{Mau01a}{article}{
  author={Maume-Deschamps, V{\'e}ronique},
  title={Correlation decay for Markov maps on a countable state space},
  journal={Ergodic Theory Dynam. Systems},
  volume={21},
  date={2001},
  number={1},
  pages={165\ndash 196},
}
\bib{Mau01b}{article}{
  author={Maume-Deschamps, V{\'e}ronique},
  title={Projective metrics and mixing properties on towers},
  journal={Trans. Amer. Math. Soc.},
  volume={353},
  date={2001},
  number={8},
  pages={3371\ndash 3389 (electronic)},
}
\bib{MelStr93}{book}{
  author={de Melo, Welington},
  author={van Strien, Sebastian},
  title={One-dimensional dynamics},
  series={Ergebnisse der Mathematik und ihrer Grenzgebiete (3) [Results in Mathematics and Related Areas (3)]},
  volume={25},
  publisher={Springer-Verlag},
  place={Berlin},
  date={1993},
  pages={xiv+605},
}
\bib{Mis81}{article}{
  author={Misiurewicz, Micha{\l }},
  title={Absolutely continuous measures for certain maps of an interval},
  journal={Inst. Hautes \'Etudes Sci. Publ. Math.},
  number={53},
  date={1981},
  pages={17\ndash 51},
}
\bib{Now88}{article}{
  author={Nowicki, Tomasz},
  title={A positive Liapunov exponent for the critical value of an $S$-unimodal mapping implies uniform hyperbolicity},
  journal={Ergodic Theory Dynam. Systems},
  volume={8},
  date={1988},
  number={3},
  pages={425\ndash 435},
}
\bib{NowSan98}{article}{
  author={Nowicki, Tomasz},
  author={Sands, Duncan},
  title={Non-uniform hyperbolicity and universal bounds for $S$-unimodal maps},
  journal={Invent. Math.},
  volume={132},
  date={1998},
  number={3},
  pages={633\ndash 680},
}
\bib{NowStr88}{article}{
  author={Nowicki, T.},
  author={van Strien, S.},
  title={Absolutely continuous invariant measures for $C\sp 2$ 
  unimodal maps satisfying the Collet-Eckmann conditions},
  journal={Invent. Math.},
  volume={93},
  date={1988},
  number={3},
  pages={619\ndash 635},
}
\bib{NowStr91}{article}{
  author={Nowicki, Tomasz},
  author={van Strien, Sebastian},
  title={Invariant measures exist under a summability condition for unimodal maps},
  journal={Invent. Math.},
  volume={105},
  date={1991},
  number={1},
  pages={123\ndash 136},
  issn={0020-9910},
}
\bib{Par60}{article}{
  author={Parry, W.},
  title={On the $\beta $-expansions of real numbers},
  language={English, with Russian summary},
  journal={Acta Math. Acad. Sci. Hungar.},
  volume={11},
  date={1960},
  pages={401\ndash 416},
}
\bib{Per74}{article}{
  author={Perekrest, V. T.},
  title={Exponential mixing in $C$-systems},
  language={Russian},
  journal={Uspehi Mat. Nauk},
  volume={29},
  date={1974},
  number={1 (175)},
  pages={181\ndash 182},
}
\bib{Pia80}{article}{
  author={Pianigiani, Giulio},
  title={First return map and invariant measures},
  journal={Israel J. Math.},
  volume={35},
  date={1980},
  number={1-2},
  pages={32\ndash 48},
  issn={0021-2172},
}
\bib{PolYur01a}{article}{
  author={Pollicott, M.},
  author={Yuri, M.},
  title={Statistical properties of maps with indifferent periodic points},
  journal={Comm. Math. Phys.},
  volume={217},
  date={2001},
  number={3},
  pages={503\ndash 520},
  issn={0010-3616},
}
\bib{Ren57}{article}{
  author={R{\'e}nyi, A.},
  title={Representations for real numbers and their ergodic properties},
  journal={Acta Math. Acad. Sci. Hungar},
  volume={8},
  date={1957},
  pages={477\ndash 493},
}
\bib{Roh48}{article}{
  author={Rohlin, V.},
  title={A ``general'' measure-preserving transformation is not mixing},
  language={Russian},
  journal={Doklady Akad. Nauk SSSR (N.S.)},
  volume={60},
  date={1948},
  pages={349\ndash 351},
}
\bib{Ros56}{article}{
  author={Rosenblatt, M.},
  title={A central limit theorem and a strong mixing condition},
  journal={Proc. Nat. Acad. Sci. U. S. A.},
  volume={42},
  date={1956},
  pages={43\ndash 47},
}
\bib{Rue68}{article}{
  author={Ruelle, D.},
  title={Statistical mechanics of a one-dimensional lattice gas},
  journal={Comm. Math. Phys.},
  volume={9},
  date={1968},
  pages={267\ndash 278},
}
\bib{Rue76}{article}{
  author={Ruelle, David},
  title={A measure associated with axiom-A attractors},
  journal={Amer. J. Math.},
  volume={98},
  date={1976},
  number={3},
  pages={619\ndash 654},
}
\bib{Rue77}{article}{ 
author={Ruelle, D.}, 
title={Applications conservant
une mesure absolument continue par rapport \`a $dx$ sur $[0,1]$},
language={French, with English summary}, 
journal={Comm.  Math.  Phys.},
volume={55}, 
date={1977}, 
number={1}, 
pages={47\ndash 51}, 
}

\bib{Ryc83}{article}{
  author={Rychlik, Marek},
  title={Bounded variation and invariant measures},
  journal={Studia Math.},
  volume={76},
  date={1983},
  number={1},
  pages={69\ndash 80},
}
\bib{Ryc88}{article}{
  author={Rychlik, Marek Ryszard},
  title={Another proof of Jakobson's theorem and related results},
  journal={Ergodic Theory Dynam. Systems},
  volume={8},
  date={1988},
  number={1},
  pages={93\ndash 109},
}
\bib{Sar02}{article}{
  author={Sarig, Omri},
  title={Subexponential decay of correlations},
  journal={Invent. Math.},
  year={2002},
  volume={150},
  number={3},
  pages={629 \ndash 653},
}
\bib{Sau00}{article}{
  author={Saussol, Beno{\^{\i }}t},
  title={Absolutely continuous invariant measures for multidimensional expanding maps},
  journal={Israel J. Math.},
  volume={116},
  date={2000},
  pages={223\ndash 248},
  issn={0021-2172},
}
\bib{Sin72}{article}{
  author={Sina\u {\i }, Ja. G.},
  title={Gibbs measures in ergodic theory},
  language={Russian},
  journal={Uspehi Mat. Nauk},
  volume={27},
  date={1972},
  number={4(166)},
  pages={21\ndash 64},
}
\bib{Tha83}{article}{
  author={Thaler, Maximilian},
  title={Transformations on $[0,\,1]$ with infinite invariant measures},
  journal={Israel J. Math.},
  volume={46},
  date={1983},
  number={1-2},
  pages={67\ndash 96},
}
\bib{ThiTreYou94}{article}{
  author={Thieullen, Ph.},
  author={Tresser, C.},
  author={Young, L.-S.},
  title={Positive Lyapunov exponent for generic one-parameter families of unimodal maps},
  journal={J. Anal. Math.},
  volume={64},
  date={1994},
  pages={121\ndash 172},
}
\bib{Tsu00a}{article}{
  author={Tsujii, Masato},
  title={Piecewise expanding maps on the plane with singular ergodic properties},
  journal={Ergodic Theory Dynam. Systems},
  volume={20},
  date={2000},
  number={6},
  pages={1851\ndash 1857},
}
\bib{Tsu01b}{article}{
  author={Tsujii, Masato},
  title={Absolutely continuous invariant measures for expanding piecewise linear maps},
  journal={Invent. Math.},
  volume={143},
  date={2001},
  number={2},
  pages={349\ndash 373},
}
\bib{Tsu93a}{article}{
  author={Tsujii, Masato},
  title={A proof of Benedicks-Carleson-Jacobson theorem},
  journal={Tokyo J. Math.},
  volume={16},
  date={1993},
  number={2},
  pages={295\ndash 310},
}
\bib{Tsu93}{article}{
  author={Tsujii, Masato},
  title={Positive Lyapunov exponents in families of one-dimensional dynamical systems},
  journal={Invent. Math.},
  volume={111},
  date={1993},
  number={1},
  pages={113\ndash 137},
  issn={0020-9910},
}
\bib{UlaNeu47}{article}{
  author={Ulam, S.},
  author={von Neumann, J.},
  title={On combination of stochastic and deterministic processes},
  date={1947},
  journal={Bull. AMS},
  volume={53},
  pages={1120},
}
\bib{Via97}{article}{
  author={Viana, Marcelo},
  title={Multidimensional nonhyperbolic attractors},
  journal={Inst. Hautes \'Etudes Sci. Publ. Math.},
  number={85},
  date={1997},
  pages={63\ndash 96},
}
\bib{Via}{book}{
  title={Stochastic dynamics of deterministic systems},
  author={Viana, Marcelo},
  series={Lecture Notes XXI Braz. Math. Colloq.},
  publisher={IMPA},
  address={Rio de Janeiro},
  year={1997},
}
\bib{Wat70}{article}{
  author={Waterman, Michael S.},
  title={Some ergodic properties of multi-dimensional $f$-expansions},
  journal={Z. Wahrscheinlichkeitstheorie und Verw. Gebiete},
  volume={16},
  date={1970},
  pages={77\ndash 103},
}
\bib{You92}{article}{
  author={Young, L.-S.},
  title={Decay of correlations for certain quadratic maps},
  journal={Comm. Math. Phys.},
  volume={146},
  date={1992},
  number={1},
  pages={123\ndash 138},
}
\bib{You98}{article}{
  author={Young, Lai-Sang},
  title={Statistical properties of dynamical systems with some hyperbolicity},
  journal={Ann. of Math. (2)},
  volume={147},
  date={1998},
  number={3},
  pages={585\ndash 650},
}
\bib{You99}{article}{
  author={Young, Lai-Sang},
  title={Recurrence times and rates of mixing},
  journal={Israel J. Math.},
  volume={110},
  date={1999},
  pages={153\ndash 188},
}

\end{biblist}
\end{bibsection}
\end{document}